\documentclass[12pt]{amsart}
\usepackage{amssymb}
\usepackage{amsthm}
\usepackage{euscript}

\theoremstyle{plain}

\newtheorem{thm}[subsection]{Theorem}

\newtheorem{cor}[subsection]{Corollary}

\begin{document}
\title[A Theorem of Kaplansky]{A theorem of Kaplansky revisited}
\author{Heydar Radjavi and Bamdad R. Yahaghi}

\address{Department of Pure Mathematics, University of Waterloo, Waterloo, Ontario, Canada N2L 3G1 \newline
\indent Department of Mathematics, Faculty of Sciences, Golestan University, Gorgan 19395-5746, Iran}
\email{hradjavi@uwaterloo.ca, bamdad5@hotmail.com,\newline  bamdad@bamdadyahaghi.com}
%\thanks{}


\bibliographystyle{plain}

\begin{abstract}
We present a new and simple proof of a theorem due to Kaplansky which unifies theorems of Kolchin and Levitzki on triangularizability of semigroups of matrices. We also give two different extensions of the theorem. As a consequence, we prove the counterpart of Kolchin's Theorem for finite groups of unipotent matrices over division rings.  We also show that the counterpart of Kolchin's Theorem over division rings of characteristic zero implies that of Kaplansky's Theorem over such division rings.

\end{abstract}

\maketitle 
%\vspace{.5cm}

\bigskip

\begin{section}
{\bf Introduction}
\end{section}

\bigskip

The purpose of this short note is three-fold. We give a simple proof of Kaplansky's Theorem on the (simultaneous) triangularizability of semigroups whose members all have singleton spectra. Our proof, although not independent of Kaplansky's,  avoids the deeper group theoretical aspects present in it and in other existing proofs we are aware of. Also, this proof can be adjusted to give an affirmative answer to Kolchin's Problem for finite groups of unipotent matrices over division rings.   In particular, it follows from this proof that the counterpart of Kolchin's Theorem over division rings of characteristic zero implies that of Kaplansky's Theorem over such division rings.  We also present extensions of Kaplansky's Theorem in two different directions. 

Let us fix some notation. Let $\Delta$ be a division ring and $M_n(\Delta)$ the algebra of all $n \times n$ matrices over $\Delta$. The division ring  $\Delta$  could    in particular be a field. By a semigroup $\mathcal{S} \subseteq M_n(\Delta)$, we mean a set of matrices closed under multiplication. An ideal $\mathcal{J}$ of $\mathcal{S}$ is defined to be a subset of  $\mathcal{S}$ with the property that $SJ \in \mathcal{J}$  
and $JS \in \mathcal{J}$ for all $ S \in \mathcal{S}$ and $ J \in \mathcal{J}$. We view the members of $M_n(\Delta)$ as linear transformations acting on the left of $\Delta^n$, where $\Delta^n$ is the right vector space of all $n\times 1$ column vectors.  A semigroup $\mathcal{S}$ is called irreducible if the orbit of any nonzero $ x \in D^n$ under $ \mathcal{S}$ spans $\Delta^n$. When $n > 1$, this is equivalent to the members of $\mathcal{S}$, viewed as linear transformations on $\Delta^n$, having no common invariant subspace other than the trivial subspaces, namely, $\{0\}$ and $\Delta^n$. On the opposite of  irreducibility is triangularizability, when the common invariant subspaces of the members of $\mathcal{S}$ include a maximal subspace chain in $\Delta^n$, i.e., there are subspaces 
$$\{0\} = \mathcal{V}_0 \subseteq  \mathcal{V}_1 \subseteq  \cdots \subset \mathcal{V}_n = \Delta^n,$$ 
where $ \mathcal{V}_j$ is a $j$-dimensional subspace invariant under every $S \in \mathcal{S}$. 

It is known, and easy to prove, that if $\mathcal{S}$ is irreducible, so is every nonzero ideal $ \mathcal{J}$ of $\mathcal{S}$ (see \cite[Lemma 2.1.10]{RR1}). Another useful fact we shall use is the continuity of spectrum for linear operators on $\mathbb{C}^n$ or $\mathbb{R}^n$; in particular, if $A$ is in the norm limit of a sequence $(A_k)_{k=1}^\infty$ of operators, then every $\lambda$ in the spectrum, $\sigma(A)$, of $A$ is the limit of a sequence $(\lambda_k)_{k=1}^\infty$ with $ \lambda_k \in \sigma(A_k)$ for each $k$ (see \cite[Lemma 3.1.2 ]{RR1}). 

Kaplansky's Theorem (see \cite[Theorem H on p. 137]{K} or \cite[Corollary 4.1.7]{RR1}) unifies two previous results: that of Levitzki, stating that a semigroup of nilpotent matrices is triangularizable (see \cite[Thoerem 35 on p. 135]{K} or \cite{L}, or \cite[Theorem 1.3]{Y2} for  a simple proof), and that of Kolchin deducing the same conclusion for a semigroup of unipotent matrices, i.e., those of the form $I+N$, where $I$ is the identity matrix and $N$ is nilpotent (see \cite{Ko} or \cite[Theorem C on p. 100]{K}). For a detailed account of Kolchin's Theorem and Kolchin's Problem, see Appendix 15B and the exercises pertaining to it on pages 60 and 171 of \cite{Ro}.

\bigskip

\begin{section}
{\bf Main Results}
\end{section}

\bigskip

\begin{thm} \label{2.1} {\bf (Kaplansky)}
Let $ n > 1$ and let $F$ be a field and $ \mathcal S$ a semigroup in $ M_n(F)$ consisting of matrices with singleton spectra. Then the semigroup 
$ \mathcal S$ is triangularizable. 
\end{thm}

\bigskip

\noindent {\bf Proof.} Without loss of generality, say by \cite[Corollary 1.3]{Y1}, we may assume that the field $F$ is algebraically closed. Now there are two cases to consider. 

(i) ${\rm ch}(F) = 0$.

We give two proofs for the assertion in this case. 

First proof. In view of Lemma 3.1.1 of \cite{RR1}, we may assume that $F = \mathbb{C}$, the field of complex numbers. If necessary, by passing to the closure of the semigroup $\mathbb{C} \mathcal{S}$,  by the continuity of spectrum, we may assume that $\mathcal{S}$ is a closed semigroup which is closed under scalar multiplications by complex numbers. Note that $ 0 \in \mathcal{S}$ and that the set of nilpotent elements of $ \mathcal{S}$ forms a semigroup ideal of $ \mathcal{S}$.   Recall that a semigroup of matrices or linear transformations is reducible iff a nonzero semigroup ideal of it is reducible. Also recall that, by Levitzki's Theorem,  every semigroup of nilpotent matrices is triangularizable, and hence reducible. Now, by passing to quotients, we only need to show that $ \mathcal S$ is reducible.  Reducibility of the semigroup is proved as soon as we show that the semigroup ideal  consisting of the nilpotent elements of $\mathcal{S}$ is nonzero. To this end, if $\mathcal{S}$ consists of scalars, the assertion trivially holds. If not, as $\mathcal{S}$ is closed under scalar multiplication by complex numbers, choose $ I + N \in \mathcal{S}$ with $ N^{k} \not= 0$ but $N^{k+1} = 0$, where $ 1 \leq k \leq n -1$. It is plain that 
$$\frac{(I+N)^n}{{n\choose k}}= \frac{1}{{n \choose k}}I + \frac{{n \choose 1}}{{n \choose k}}N + \cdots + \frac{{n \choose k-1}}{{n \choose k}}N^{k-1} + N^k, $$ 
for all $ n \in \mathbb{N}$  with $ n > k$. This clearly gives 
$$ N^k = \lim_n \frac{(I+N)^n}{{n\choose k}} \in \mathcal{S}, $$
showing that the semigroup ideal of nilpotent elements of $\mathcal{S}$ is nonzero, as desired. 

Second proof. By passing to $F^*\mathcal{S}$, where $F^* = F \setminus \{0\}$, we may assume that $\mathcal{S}$ is closed under scalar multiplications by the nonzero elements of $F$. Again, we only need to show that $\mathcal S$ is reducible. If the semigroup $\mathcal{S}$ contains a nilpotent element, then reducibility of $\mathcal{S}$ follows from that of the nonzero semigroup ideal of  all nilpotent elements of  $\mathcal{S}$. So it remains to prove the assertion when $\mathcal{S}$ contains no nonzero nilpotent element.  It is then plain that $\mathcal{S}$ is reducible iff the set of all unipotent elements of $\mathcal{S}$ is reducible. Thus, in view of Kolchin's Theorem, we will be done as soon as we prove that the set of all unipotent elements of $\mathcal{S}$ forms a semigroup. To this end, let $ I + N_1 , I + N_2 \in \mathcal{S}$ be arbitrary unipotent elements. We can write $ (I+ N_1)(I + N_2) = cI + N' \in \mathcal{S}$, where $c \in   F^*$ and  $N' $ is a nilpotent matrix. We need to show that $c=1$. If $N_2 = 0$, we have nothing to prove. Suppose $N_2 \not= 0$ so that $N_2^k \not= 0$ but $N_2^{k+1} = 0$ for some $k \in \mathbb{N}$ with $k < n$. Thus
$$(I + N_2)^m = I + {m \choose 1} N_2 + \cdots + {m \choose k} N_2^k, $$
for all $ m \in \mathbb{N}$. Recall that ${m \choose k} := 0$ whenever $ m <k$. Clearly, we have $ (I+ N_1)(I + N_2)^m = c_mI + N_m' \in \mathcal{S}$  for all $ m \in \mathbb{N}$ with $c_m \in F$ an $n$th root of unity, i.e., $c_m^n = 1$ and $N'_m$ a nilpotent matrix. Since the set of the $n$th roots of unity in $F$ has at most  $n$ elements, we see that there exists a subsequence $(m_i)_{i=1}^\infty$  and an $l \in \mathbb{N}$ such that $ (I+ N_1)(I + N_2)^{m_i} = c_{m_l}I + N'_{m_l}$ for all $ i \in \mathbb{N}$.
Therefore, 
$$ \Big( (I+ N_1)(I + N_2)^{m} - c_{m_l}I\Big)^n = 0, $$
for infinitely many $ m \in \mathbb{N}$. Now, fix $ 1 \leq i , j \leq n$ and note that the $(i, j)$ entry of the matrix $ \Big( (I+ N_1)(I + N_2)^{m} - c_{m_l}I\Big)^n$ is a polynomial of degree $k$ in $m$ having infinitely many roots, namely, $m_j$'s ($j \in \mathbb{N}$). This implies that   the $(i, j)$ entry of the matrix $ \Big( (I+ N_1)(I + N_2)^{m} - c_{m_l}I\Big)^n$ is zero for all $m \in \mathbb{N} \cup \{0\}$. Consequently, 
$$ \Big( (I+ N_1)(I + N_2)^{m} - c_{m_l}I\Big)^n = 0, $$
for all $ m \in \mathbb{N} \cup \{0\}$. Setting $m=0 ,1$ in the above, we obtain $ c = c_{m_l}=1$, which is what we want.

(ii) ${\rm ch}(F) \not= 0$. 

Let $ p = {\rm ch}(F)$. Again, if necessary, by passing to $F^*\mathcal{S}$, where $F^* = F \setminus \{0\}$, we may assume that $\mathcal{S}$ is closed under scalar multiplications by the nonzero elements of $F$. Again, we only need to show that $\mathcal S$ is reducible. To this end, as the set of nilpotent elements of $ \mathcal{S}$ forms a semigroup ideal of $ \mathcal{S}$, the assertion follows from Levitzki's Theorem as soon as $\mathcal{S}$ contains a nonzero nilpotent. However, if $\mathcal{S}$ does not contain any nonzero nilpotent,   set 
$$\mathcal{S}_1 := \{ S \in \mathcal{S} : \det S = 1 \}.$$
Clearly, $\mathcal{S}$ is reducible if and only if  $\mathcal{S}_1$  is reducible. By way of contradiction, suppose on the contrary that $\mathcal{S}_1$ is irreducible.  Note that  $\mathcal{S}_1$ is a subsemigroup of $\mathcal{S}$. In fact $\mathcal{S}_1$ forms a group of matrices. To see this, let $ c I + N \in \mathcal{S}_1$ with $c^n=1$ be arbitrary. Choose $ r \in \mathbb{N}$ such that $ p^r > n$ and note that 
$ (cI + N)^{p^r} = c^{p^r} I$. This clearly yields 
$$ c^{-1} (cI + N)^{-1} = \left(I + \frac{N}{c}\right)^{-1}=  \left(I + \frac{N}{c}\right)^{p^r-1} \in \mathcal{S},$$
which implies $ (cI + N)^{-1} \in \mathcal{S}$. But 
$$ \det (cI + N)^{-1} = \Big(\det(cI + N)\Big)^{-1} = c^{-n} = 1.$$ 
Thus $ (c I + N)^{-1} \in \mathcal{S}_1$. Now, let $ S= cI + N  \in \mathcal{S}_1$ be arbitrary. It follows that $ \det S = c^n = 1$, and hence $ c \in \Omega$, where $\Omega = \{ \omega_1 , \ldots, \omega_n\}$ denotes the set of $n$th roots of  unity in the algebraically closed field $F$. It thus follows that $ {\rm tr} (S) \in n \Omega$ for all $ S \in \mathcal{S}_1$. Since $n \Omega$ is a finite set and $\mathcal{S}_1$ is irreducible, it follows from Theorem B on page 99 of \cite{K} (or from \cite[Theorem 1]{RR2}) that $\mathcal{S}_1$ is a finite group. We now see that $ \frac{\mathcal{S}_1}{Z(\mathcal{S}_1)}$, where $ Z(\mathcal{S}_1) = FI \cap \mathcal{S}_1$ denotes the center of $\mathcal{S}_1$, is a finite $p$-group. Note that the center of $\mathcal{S}_1$ is a subset of  $FI$ because $\mathcal{S}_1$ is irreducible. It is a well-known fact that the center of any finite $p$-group is nontrivial. Thus, there exists an $ A = c I + N \in  \mathcal{S}_1 \setminus FI $ with the property that for all $ B \in \mathcal{S}_1$ we can find a $c_B \in F$ such that $ AB = c_B BA$. Let $ \mathcal{N} = \ker (A - c I)$. Clearly, $ \mathcal{N}$ is a nontrivial subspace of $F^n$. We obtain a contradiction by showing that $ \mathcal{N}$ remains invariant under  $\mathcal{S}_1$. To this end, let  $ B \in \mathcal{S}_1 $ and  $ x \in \mathcal{N}$ be arbitrary. If $Bx = 0$, then $ Bx \in \mathcal{N}$. If $Bx \not= 0$, then we can write 
$$ AB x = c_B BA x = c_B c Bx, $$ 
implying that  $ c_B c $ is an eigenvalue for $A$. But $c$ is the only eigenvalue of $A$. This means $ c_B c = c$, implying that $c_B = 1$, for $ c \not= 0$. Consequently, $ABx = cBx$, which means $ Bx \in  \mathcal{N}$. That is, $ \mathcal{N}$ is a nontrivial invariant subspace for $\mathcal{S}_1$, a contradiction. This completes the proof. 
\hfill \qed

\bigskip 

\noindent {\bf Remark.} In fact, in view of the Gordon-Motzkin Theorem (\cite[Theorem 16.4]{La}), our second proof of Kaplansky's Theorem in the case of ${\rm ch}(F) = 0$ shows that the counterpart of Kolchin's Theorem over division rings of characteristic zero implies that of Kaplansky's Theorem over such division rings. In other words, if every semigroup of unipotent matrices over a division ring  $\Delta$ of characteristic zero is triangularizable, then so is every semigroup of matrices of the form $ cI + N$, where $c$ comes from the center of $\Delta$ and $ N$ is a nilpotent matrix with entries from $\Delta$. 

\bigskip

In fact the proof of the theorem establishes the following counterpart of Kaplansky's Theorem for finite semigroups of matrices over division rings. In particular, the following proves the counterpart of Kolchin's Theorem for finite groups of unipotent matrices over division rings (see page 62 of \cite{Ro}). It is worth mentioning that the corresponding statement is known to be true for arbitrary groups of unipotent matrices over division rings if the characteristic is zero or large enough (see \cite{M}).

\bigskip

\begin{thm} \label{2.2}
Let $ n > 1$ and let $\Delta$ be a division ring  and $ \mathcal S$ a finite semigroup in $ M_n(\Delta)$ consisting of matrices of the form $ cI + N$, where $ c$ is in the center of $\Delta$ and $N$ is nilpotent. Then the semigroup  $ \mathcal S$ is triangularizable. 
\end{thm}

\bigskip

\noindent {\bf Proof.} By passing to quotients, we only need to prove reducibility. If  
$ \mathcal S$  contains a nonzero nilpotent matrix, the reducibility of $ \mathcal S$ follows from that of its nonzero semigroup ideal of its nilpotent elements. So we may without loss of generality assume that the finite semigroup $ \mathcal S$ consists of invertible matrices. This clearly implies that $ \mathcal S$ is indeed a group of matrices. Now the assertion is easy if the characteristic of $\Delta$ is zero. That is because every element $cI + N$ of the finite group $\mathcal S$ has finite order. This easily implies that $N=0$, from which reducibility of $\mathcal S$ follows. Next suppose $\Delta$ has a nonzero characteristic $p$. We now see that $ \frac{\mathcal{S}}{Z(\mathcal{S})}$, where $ Z(\mathcal{S}) \subset FI \cap \mathcal{S}$ denotes the center of $\mathcal{S}$, is a finite $p$-group. That is because $ (cI + N)^{p^r} = c^{p^r} I$ for all $A= cI + N \in \mathcal{S}$, where $ r \in \mathbb{N}$ is chosen such that $ p^r > n$. 
Note that the center of $\mathcal{S}$ is a subset of  $FI$ because $\mathcal{S}$ is irreducible. Again the center of $\mathcal{S}$ is nontrivial because it is a finite $p$-group. Thus, there exists an $ A = c I + N \in  \mathcal{S} \setminus FI $ with the property that for all $ B \in \mathcal{S}$ a $c_B \in F$ can be found such that $ AB = c_B BA$. Let $ \mathcal{N} = \ker (A - c I)$. Clearly, $ \mathcal{N}$ is a nontrivial subspace of $\Delta^n$. We obtain a contradiction by showing that $ \mathcal{N}$ remains invariant under  $\mathcal{S}$. To this end, let  $ B \in \mathcal{S} $ and  $ x \in \mathcal{N}$ be arbitrary. If $Bx = 0$, then $ Bx \in \mathcal{N}$. If $Bx \not= 0$, then we can write 
$$ AB x = c_B BA x = c_B c Bx, $$ 
implying that  $ \Big((c - c_B c)I + N \Big)Bx =0 $, which in turn implies that $(c - c_B c)I + N$ is not invertible. This yields $ c = c_B c$ because $N$ is nilpotent. Thus $c_B = 1$, for $ c \not= 0$. Consequently, $ABx = cBx$, which means $ Bx \in  \mathcal{N}$. That is, $ \mathcal{N}$ is a nontrivial invariant subspace for $\mathcal{S}_1$, a contradiction. This completes the proof.
\hfill \qed

\bigskip 

Theorems \ref {2.3} and \ref{2.5} below are extensions of Kaplansky's Theorem. 

\bigskip

\begin{thm} \label{2.3}
Let $ n > 1$  and let $ \mathcal S$ be a semigroup in $ M_n(\mathbb{C})$ with the property that $ \sigma (S) \subset r_S \mathbb{T}$ for all $S \in \mathcal{S}$, where $ r_S \geq 0$ and $\mathbb{T}$ denotes the unit circle in the complex plane. Then the semigroup  $ \mathcal S$ is simultaneously similar to a block triangular semigroup in which each diagonal block of size greater than $1$ is an irreducible semigroup of the form $ \mathbb{R}_{\geq 0} \mathcal{S}_u$, where $ \mathbb{R}_{\geq 0} = [0, +\infty)$ and $\mathcal{S}_u$ is an irreducible semigroup of unitaries. 
\end{thm}

\bigskip

\noindent {\bf Proof.} If necessary, replacing $ \mathcal{S}$ by $ \mathbb{R}^+ \mathcal{S}$, we may assume that $ \mathcal{S}= \mathbb{R}^+ \mathcal{S}$. Applying a simultaneous similarity, we can assume that the number $k$ of the diagonal blocks is maximal so that each diagonal block is irreducible. Let $ \mathcal{S}_0$ denote an arbitrary diagonal block of size $n_0 >1$. Clearly $ \mathcal{S}_0$ is  an irreducible semigroup in $ M_{n_0}(\mathbb{C})$ and for all 
$S \in \mathcal{S}_0$, there is an $ r_S \geq 0$  such that $ \sigma (S) \subset r_S \mathbb{T}$. If there is a nonzero nilpotent in $ N \in \mathcal{S}_0$, then $ \mathcal{S}_0 N \mathcal{S}_0$ is a nonzero semigroup ideal of nilpotents, which is triangularizable by Levitzki's Theorem, and hence reducible, a contradiction. Thus, $\mathcal{S}_0$ contains no nonzero nilpotents. Hence $ \mathcal{S}_0 \setminus \{0\} $
is a semigroup of invertibles.  Now let 
$$\mathcal{S}_1 = \{ S \in \mathcal{S}_0:   |\det S| = 1\}.$$ 
It is now clear that $\mathcal{S}_1$ is irreducible, that $ \sigma (S) \subseteq \mathbb{T}$ for all $S \in \mathcal{S}_1$, and that  $\mathcal{S}_0 \subseteq \mathbb{R}_{\geq 0}  \mathcal{S}_1$. Consequently, the spectrum, and hence the trace functional, is bounded on $\mathcal{S}_1$. It thus follows from \cite[Theorem 4]{RR2} that the semigroup $\mathcal{S}_1$ is bounded. Now, in view of \cite[Theorem 3.1.5]{RR1}, if necessary, by passing to the closure of  $\mathcal{S}_1$, which is a bounded group of matrices,  we see that  $\mathcal{S}_1$  is simultaneously similar to a semigroup of unitaries. This completes the proof  because $ \mathcal{S}_0$ was an arbitrary diagonal block of size greater than $1$.
\hfill \qed

\bigskip

\begin{cor} \label{2.4}
Let $ n > 1$  and let $ \mathcal U$ be a group of unitary matrices in $ M_n(\mathbb{C})$. Let $\mathcal{S}$ be a semigroup each of whose members is of the form $cU + N$, where $ U \in \mathcal{U}$, $c \in \mathbb{C}$, and $N$ is a nilpotent matrix such that $\{ U, N\}$ is triangularizable.  Then the conclusion of Theorem \ref{2.3} holds. 
\end{cor}

\bigskip

\noindent {\bf Proof.} Note that for every member $S = c_SU_S + N_S \in \mathcal{S}$ with $ U_S \in \mathcal{U}$, $c_S \in \mathbb{C}$, and $N_S$ nilpotent, we have $\{ U_S, N_S\}$ is triangularizable. Thus, the semigroup  $\mathcal{S}$ has the property that $ \sigma (S) \subset r_S \mathbb{T}$ for all $S \in \mathcal{S}$, where $ r_S= |c_S| \geq 0$ and $\mathbb{T}$ denotes the unit circle in the complex plane. Therefore, the assertion follows from the preceding theorem. 
\hfill \qed

\bigskip 

\noindent {\bf Remark.} It is now clear why Theorem \ref{2.3} is an extension of Kaplansky's Theorem. Let  $ \mathcal S$ be a Kaplansky semigroup in $M_n(F)$, namely a semigroup in $ M_n(F)$ consisting of matrices of the form $ cI + N$, where $ c \in F$  and $N$ is nilpotent, or equivalently $\mathcal{S}$ is a semigroup with $\mathcal{U} = \{I\}$ in the preceding corollary. We only need to show that $\mathcal{S}$ is reducible. If not, then by the preceding corollary $\mathcal{S}$ is similar to an irreducible semigroup of the form $ \mathbb{R}^+ \mathcal{S}_u$, where $\mathcal{S}_u$ is an irreducible semigroup of unitary matrices. That is, we have $S = r_S U_S$ for all $S \in \mathcal{S} \setminus \{0\}$ with $r_S = |c_S|> 0$ and $U_S$ unitary. Since $\sigma (S) = \{c_S\}$, we see that $U_S$ is a scalar matrix for all $S \in \mathcal{S} \setminus \{0\}$  because it is diagonalizable and has singleton spectrum. Consequently, $\mathcal{S}$ is reducible, contradicting the contradiction hypothesis. This proves Kaplansly's Theorem.

\bigskip

For a collection $\mathcal{C}$ in $M_n(F)$, by the commutant of $\mathcal{C}$, denoted by $\mathcal{C}'$, we mean 
$$\mathcal{C}' := \{A \in M_n(F) : AB = BA \   \forall \  B \in \mathcal{C}'\}.$$

\bigskip

\begin{thm} \label{2.5}
Let $ n \in \mathbb{N}$ and let $F$ be a field, $\mathcal{A}$ a commutating set of triangularizable  matrices in $ M_n(F)$, $\mathcal{N}$ the set of all nilpotents in $ M_n(F)$, and $ \mathcal S$ a semigroup in $ M_n(F)$ consisting of matrices of the form $A + N$, where $ A \in \mathcal{A}$ and $ N \in \mathcal{A}' \cap \mathcal{N}$. Then the semigroup  $ \mathcal S$ is triangularizable. 
\end{thm}

\bigskip

\noindent {\bf Proof.}  We view the elements of $M_n(F)$ as linear transformations on $F^n$ and proceed by induction on $n$, the dimension of the underlying space. The assertion trivially holds for $n=1$. Assume $ n > 1$ and that the assertion holds for such semigroups of linear transformations acting on spaces of dimension less than $n$. If $\mathcal{A}$ consists of matrices with sinlgeton spectra, then every $S \in \mathcal{S}$ is of the form $S= A+ N $, where $ A= cI + N_1$ and $ N \in \mathcal{A}' \cap \mathcal{N}$. Thus, $N$ commutes with $A$, and hence with $N_1$, and thus $N_1 + N$ is nilpotent as well. Consequently, every $S \in \mathcal{S}$ has singleton spectrum. So the assertion follows from Kaplansky's Theorem in this case. If not, choose $S_0 = B + M \in \mathcal{S}$, where $B \in \mathcal{A}$ has more than one point in its spectrum and $ M \in \mathcal{A}' \cap \mathcal{N}$. Therefore, there exists a nontrivial direct sum decomposition  $F^n = \mathcal{V}_1 \oplus \mathcal{V}_2$ with respect to which $B= B_1 \oplus B_2$ with $\sigma(B_1) \cap \sigma(B_2) = \emptyset$ with $B_1 = B|_{\mathcal{V}_1}$ and $B_2 = B|_{\mathcal{V}_2}$. Note that $ \mathcal{A} \cup (\mathcal{N} \cap \mathcal{A}') \subseteq \{B\}'$. For each $S \in \mathcal{S}$,  with respect to the decomposition $F^n = \mathcal{V}_1 \oplus \mathcal{V}_2$, we can write $S= A_S + N_S = (A_1 + N_1) \oplus (A_2 + N_2)$ where $A_j, N_j \in \{B_i\}' $ ($j =1, 2$). 
For $j = 1, 2$, let 
$$\mathcal{S}_j = \mathcal{S}|_{\mathcal{V}_j} = \{ A_S|_{\mathcal{V}_j} + N_S|_{\mathcal{V}_j} : S \in \mathcal{S} \}$$
and let $ \mathcal{N}_j$ denote the set of all nilpotent linear transformations on  $\mathcal{V}_j$. 
Note that for $ j = 1, 2$, $\mathcal{S}_j$ is a semigroup of linear transformations  on $\mathcal{V}_j$ and  $\mathcal{A}|_{\mathcal{V}_j}$ is commutative and consists of triangularizable linear transformations.  
%$ \mathcal{N}_j \cap  \mathcal{A}'|_{\mathcal{V}_j} \subseteq (\mathcal{A}|_{\mathcal{V}_j})'$
Let $A_j + N_j \in \mathcal{S}_j$ with $j \in \{1, 2\}$ be an arbitrary element so that $S= A_S + N_S = (A_1 + N_1) \oplus (A_2 + N_2)$ for some $S \in \mathcal{S}$. It follows that 
$N_1 \oplus  N_2 \in \mathcal{A}' \cap \mathcal{N}$, and hence
$ N_j \in (\mathcal{A}|_{\mathcal{V}_j})' \cap \mathcal{N}_j$ for $j = 1, 2$. Consequently, every element of $\mathcal{S}_j$ is of the form $A_j + N_j$, where $A_j \in \mathcal{A}|_{\mathcal{V}_j}$ and $N_j \in (\mathcal{A}|_{\mathcal{V}_j})' \cap \mathcal{N}_j$ for each $j = 1, 2$. 
  So by the induction hypothesis $\mathcal{S}_1$ and $\mathcal{S}_2$, and hence 
$\mathcal{S}$, are triangularizable. This completes the proof. 
%If necessary, applying a similarity, we may assume that $ B = \oplus_{i=1}^k B_i$, where $ B_i = \alpha_i I + M_i$, where $\alpha_i$'s ($1 \leq i \leq k$) are in $F$ and distinct, $I$ stands for the identity matrix of an appropriate size, and $M_i$ is nilpotent. Note that if $A+ N \in \mathcal{A}$, then $A , N \in \{B\}'$, where $\{B\}' := \{A \in M_n(F) : AB = BA\}$ denotes the commutant of $\{B\}$. If $A+ N \in \mathcal{A}$, as  $\alpha_i$'s ($1 \leq i \leq k$)  and pairwise distinct, we see that $A  = \oplus_{i=1}^k A_i$ and $N  = \oplus_{i=1}^k N_i$ with $A_i, N_i \in \{B_i\}' = \{M_i\}'$.  
\hfill \qed

\bigskip

\begin{cor} \label{2.6}
Let $ n > 1$ and let $F$ be a field, $\mathcal{A}$ a diagonalizable set of matrices in $ M_n(F)$ and $ \mathcal S$ a semigroup in $ M_n(F)$ consisting of matrices of the form $A + N$, where $ A \in \mathcal{A}$ and $ N $ is in the commutant of $\mathcal{A}$ and is nilpotent. Then the semigroup  $ \mathcal S$ is triangularizable. 
\end{cor}

\bigskip

\noindent {\bf Proof.} Just note that every diagonalizable set of matrices is commutative and triangularizable. Thus the assertion is a quick consequence of Theorem  \ref{2.5}. 
\hfill \qed

\bigskip 

We conclude with the following question which we have not been able to resolve. 

\bigskip

\noindent {\bf Question.} Let $ n > 1$ and let $F$ be a field and $\mathcal{T}$ a triangularizable set of  matrices in $ M_n(F)$ and $ \mathcal S$ a semigroup in $ M_n(F)$ consisting of matrices of the form $T + N$, where $ T \in \mathcal{T}$ and $ N $ is in the commutant of $\mathcal{T}$ and is nilpotent. Is the semigroup  $ \mathcal S$ necessarily triangularizable?

\bigskip

\noindent  {\bf Acknowledgment.}  The authors would like to thank the referee for reading the paper carefully and making helpful comments. This paper was submitted while the second-named author was on sabbatical leave, graciously hosted by the Department of Pure Mathematics of the University of Waterloo.


\vspace{2cm}
\begin{thebibliography}{999}




%\bibitem[1]{B}
%{\large [B]}
%W. Burnside, On the condition of reducibility of any
%group of linear substitutions, {\it Proc. London Math. Soc.} 3 (1905), 430-434.


%\vspace{1mm}
%\bibitem[2]{D}
%{\large [DS]}
%P.K. Draxl, {\it Skew Fields}, Cambridge University Press, 1983.




%\vspace{1mm}
%\bibitem[3]{H}
%{\large [H]}
%T.W. Hungerford, {\it Algebra}, Springer Verlag, New York, 1974.


\vspace{1mm}
\bibitem[1]{K}
I. Kaplansky, {\it  Fields and Rings} 2nd ed, University of Chicago Press, Chicago, 1972.


\vspace{1mm}
\bibitem[2]{Ko}
%{\large [Ko]}
E. Kolchin, On certain concepts in the theory of algebraic matric groups,  {\it  Ann. of Math.} Vol. {\bf 49} (4), 1948, 774-789.


\vspace{1mm}
\bibitem[3]{La}
T.Y.  Lam, {\it  A first course in noncomutative rings}, Springer Verlag, New York, 1991.



\vspace{1mm}
\bibitem[4]{L}
%{\large [L]}
J. Levitzki, {\"U}bber Nilpotente Unterringe,  {\it  Math. Ann.} Vol. {\bf 105}, 1931, 620-627.





\vspace{1mm}
\bibitem[5]{M}
%{\large [RY]}
H.D. Mochizuki, Unipotent matrix groups over division rings,  {\it Canadian Math. Bull.} Vol. {\bf 21} (2), 1978, 249-250.




\vspace{1mm}
\bibitem[6]{RR1}
%{\large [RR]}
H. Radjavi and P. Rosenthal, {\it Simultaneous
Triangularization}, Springer Verlag, New York, 2000.


\vspace{1mm}
\bibitem[7]{RR2}
%{\large [RR]}
H. Radjavi and P. Rosenthal, Limitations on the size of semigroups of matrices, Semigroup Forum,  76: 25-31, 2008.




\vspace{1mm}
\bibitem[8]{Ro}
%{\large [RR]}
L.H. Rowen, {\it Graduate Algebra: Noncommutative View}, American Mathematical Society, Providence, RI, 2008.




\vspace{1mm}
\bibitem[9]{Y1}
%{\large [Y1]}
B.R. Yahaghi, On irreducible semigroups of martices with
traces in a subfield, {\it Linear Algebra Appl.} {\bf 383} (2004), 17-28.



\vspace{1mm}
\bibitem[10]{Y2}
%{\large [Y2]}
B.R. Yahaghi, On $F$-algebras of algebraic matrices over a subfield $F$ of the center
 of a division ring, {\it  Linear Algebra Appl.} {\bf 418} (2006),
599-613.

%\vspace{1mm}
%\bibitem[14]{Y3}
%{\large [Y3]}
%B.R. Yahaghi, {\it Reducibility Results on Operator
%Semigroups}, Ph.D. Thesis, Dalhousie University, Halifax, Canada, 2002.



\end{thebibliography}
\end{document}